\input amstex
\documentstyle{amsppt}
\document
\topmatter
\title
Neutral bi-Hermitian Gray surfaces.
\endtitle
\author
Wlodzimierz Jelonek
\endauthor

\abstract{The aim of this paper is to give examples of compact neutral 4-manifolds  $(M,g)$ whose Ricci tensor $\rho$  satisfies the relation $\nabla_X\rho(X,X) =\frac13X\tau g(X,X)$.  We present also a family of new Einstein bi-Hermitian neutral metrics on ruled surfaces of genus $g>1$. }
 \endabstract

\bigskip
\thanks{MS classification:53C05,53C20,53C25. The paper was partially supported by KBN grant  2P0 3A 02324.}\endthanks

\endtopmatter
\define\G{\Gamma}
\define\DE{\Cal D^{\perp}}
\define\e{\epsilon}
\define\n{\nabla}
\define\om{\omega}

\define\k{\diamondsuit}
\define\th{\theta}

\define\a{\alpha}

\define\lb{\lambda}
\define\A{\Cal A}
\define\AC{\Cal AC^{\perp}}

\define\1{\Cal D_{\lb}}
\define\2{\Cal D_{\mu}}
\define\0{\Omega}

\define\bJ{\bar J}

\define\Si{\Sigma}

\define\De{\Cal D}

{\it This paper I dedicate to Professor Kouei Sekigawa on his sixtieth birthday.    }

\vskip 3cm

\bigskip
{\bf 0. Introduction. } In the present paper we are concerned with the class of neutral semi-Riemannian 4-manifolds $(M,g)$ whose Ricci tensor $\rho$    satisfies the condition
$$\n_X\rho(X,X) =\frac13X\tau g(X,X)\tag *$$
where $\tau$ is the scalar curvature of $(M,g)$. These class of manifolds was introduced by A. Gray ( see [G],[Be]).  For general facts and some results concerning neutral 4-manifolds we refer to [K], [Pe],[A], [D],[M-L]. In [M] the condition for the existence of two opposite almost complex structures on 4-manifolds are studied.  Note that in [A] (p.187) it is stated that there are no known
neutral Hermitian-Einstein metrics on ruled surfaces except trivial ones. In this note we give among others many new, non-trivial  examples of such manifolds. The main subject of the paper is the construction of new Einstein and Gray neutral metrics on ruled surfaces. The methods however are very similar to those used by me in [J-1]-[J-4].

In our papers [J-2]-[J-4] we have described  Riemannian metrics $g$ on compact complex surfaces $(M,J)$ such that $(M,g)$ satisfies the condition $*$ and has $J$-invariant  Ricci tensor.
In particular we have given a complete classification of  bi-Hermitian Gray surfaces on ruled surfaces of genus $g\ge 1$ and surfaces which are of co-homogeneity 1 on ruled surfaces of genus $g=0$. We denote by $g$ also the Riemannian metric but it should not cause any misunderstandings. The technics used in [J-3] can be applied to describe a large class of neutral bi-Hermitian Gray surfaces.  The equations, describing such surfaces whose Ricci tensor has exactly two real eigenvalues such that the corresponding two dimensional eigendistributions are space-like and time-like, are described by the same equations as in [J-3],[J-4]. In that way we present a large class of explicit examples of neutral bi-Hermitian Gray surfaces.
In particular we get a large class of Einstein neutral bi-Hermitian  surfaces on ruled surfaces of genus $g>1$.
\bigskip
{\bf 1. Neutral $\AC$-surfaces.}
By an $\AC$- manifold (see [Be],[G]) we mean a semi-Riemannian manifold $(M,g)$ satisfying the condition
$$\frak C_{X Y Z}\n_X\rho(Y,Z)=\frac 2{(\text{dim}M+2)}\frak C_{X Y Z}X\tau g(Y,Z),  \tag *$$
where $\rho$ and $\tau$ are the Ricci tensor and the scalar curvature  of $(M,g)$ respectively and $\frak C$ means the cyclic sum. In this paper,
an $\AC$-manifold with neutral metric is also called a neutral Gray manifold.
A Riemannian manifold $(M,g)$ is an $\AC$ manifold if and only if the Ricci endomorphism $Ric$ of $(M,g)$ is of the form $Ric=S+\frac{2}{n+2}\tau Id$ where $S$ is a Killing tensor and $n=$dim$M$. Let us recall that a symmetric (1,1) tensor $S$ on a semi-Riemannian manifold $(M,g)$ is called a Killing tensor if $g(\n S(X,X),X)=0$ for all $X\in TM$ and that a semi-Riemannian manifold whose Ricci tensor is a Killing tensor is called an $\A$-manifold.

On a semi-Riemannian manifold $(M,g)$ a distribution $\De\subset TM$ is called umbilical (see [J-3]) if $\De$ is non-degenerate (i.e. a metric $g$ is non-degenerate on $\De$) and
$\n_XX_{|\DE}=g(X,X)\xi$ for every $X\in\G(\De)$, where $X_{|\DE}$ is the $\DE$ component of $X$ with respect to the orthogonal decomposition $TM=\De\oplus\DE$. The vector field $\xi$ is called the mean curvature normal of $\De$.  The foliation tangent to  involutive distribution $\De$  is called totally geodesic if its every leaf is a totally geodesic ( i.e $\n_XX\in \De$ for any  $X \in\De$ )
non-degenerate submanifold of $(M,g)$. In the sequel we shall not distinguish between $\De$ and a  foliation tangent to $\De$ and we shall also say that $\De$ is totally geodesic in such a case.

 It is not difficult to prove exactly as in the Riemannian case (see [J-3]) the following lemma:
\medskip
{\bf Lemma  1. } {\it Let $S\in End(TM)$ be a (1,1) tensor on a neutral semi-Riemannian 4-manifold $(M,g)$. Let us assume that $S$ has exactly two everywhere different eigenvalues $\lb,\mu$ of the same multiplicity 2,  i.e. \text{ dim} $\1$=\text{ dim }$\2=2$, where $\1,\2$ are non-degenerate eigendistributions of $S$ corresponding to $\lb,\mu$ respectively. Then $S$ is a Killing tensor if and only if both distributions $\1$ and $\2$ are umbilical with mean curvature normal equal respectively}
$$\xi_{\mu}=\frac{\n\mu}{2(\lb-\mu)}, \  \xi_{\lb}=\frac{\n\lb}{2(\mu-\lb)}.$$

 We shall call a bi-Hermitian surface with bi-Hermitian Ricci tensor simply as bi-Hermitian surface If $(M,g)$ is also an $\AC$ manifold, then we call it a neutral bi-Hermitian Gray surface.
\medskip
{\bf Proposition 1. } {\it Let us assume that $(M,g)$ is a  simply connected neutral Gray 4-manifold and  the Ricci tensor $S$ ($\rho(X,Y)=g(SX,Y)$) has exactly two real eigenvalues $\lb,\mu$ and no null eigenvectors.   Then there exist two Hermitian  complex structures $J,\bJ$ commuting with $S$ and $(M,g)$ is a bi-Hermitian neutral Gray surface.}
\medskip
{\it Proof. }  Analogous to [J-1].$\k$

We shall prove
\bigskip
{\bf Proposition 2.} {\it Let $(M,g)$ be a 4-dimensional neutral semi-Riemannian manifold. Let $\De$ be a two dimensional totally geodesic bundle-like and space-like foliation on $M$. Then $M$ admits (up to 4-fold covering) two opposite Hermitian structures $J,\bJ$ such that $J|{\De}=-\bJ|{\De},J|{\DE}=\bJ|{\DE}$. The nullity distributions of both $J,\bJ$ contain $\De$.}
\medskip
{\it Proof.} To prove the first part of the Proposition it is enough to show that on $M$ there exists a Killing tensor with eigendistributions $\De,\DE$. Since the foliation $\De$ is  totally geodesic and bundle-like it follows that $\n_XX\in \G(\De) \text{ (resp.) }\n_XX\in\G(\DE)$ if $X\in \G(\De) \text{ (resp.) }X \in \G(\DE)$. Consequently a tensor $S$ defined by
$$\gather
SX=\lb X\text{ if } X\in \De\\
SX=\mu X\text{ if } X\in \DE\endgather$$
where $\lb\ne\mu$ are two different real numbers is a smooth Killing tensor.  We can assume (up to 4-fold covering) that   the distributions $\De,\DE$ are orientable. Let us denote by $J$ the only almost Hermitian structure which preserves $\De,\DE$ and agrees with their orientations.
We define $\bJ$ by $J|{\De}=-\bJ|{\De},J|{\DE}=\bJ|{\DE}$. From Proposition 1 it follows that both
$J,\bJ$ are Hermitian. Now let $\{E_1,E_2\}$ be a local orthonormal frame on $\De$. Then
$$\n J(E_1,E_1)+J(\n_{E_1}E_1)=\n_{E_1}E_2.$$
It follows that $\n J(E_1,E_1)\in\De$ and consequently $\n_{E_1}J=0$. Analogously $\n_{E_2}J=0.$
$\k$
\bigskip

Now we give a theorem, whose proof is analogous as in the Riemannian case ( see [M-S],[B],[S]).

\bigskip
{\bf Theorem  1.}  {\it Let us consider the  manifold
$ U=(a,b)\times P$,
where $(P,g_0)$ is a 3-dimensional semi-Riemannian $\A$-manifold (a principle  $S^1$ bundle $p:P\rightarrow \Si$ ) over a Riemannian surface $(\Si,g_{can})$ of constant sectional curvature  $K$, with an $\AC$-metric
$$g=dt^2+f(t)^2\th^2-g(t)^2p^*g_{can},\tag 1.5$$
where $g_0=\th^2-p^*g_{can}$, $f,g\in C^{\infty}(a,b)$ and $\th$ is the connection form of $P$. Then the metric $g$ on $U$ extends to a smooth $\AC$ metric on the ruled surface $M$ which is a $\Bbb CP^1$-bundle over $\Si$ and such that $U$ is an open and dense subset of $M$ if and only if the functions $f,g\in C^{\infty}(a,b)$ satisfy the conditions:

(a) $f(a)=f(b)=0, f'(a)=1, f'(b)=-1, $

(b)   $g(a)\ne0\ne g(b),g'(a)=g'(b)=0. $ }

\medskip
{\it Remark.}  Let us note that the metric (1.5) induces a semi-Riemannian metric on $M$ if the functions $f,g$ satisfy:

(i) $f$ is positive on $(a,b)$ , and $f$ is odd at $a$ and $b$ , i.e. $f$ is the restriction of a function $f$ on $\Bbb R$ satisfying $  f(a+t)=-f(a-t), f(b+t)=-f(b-t)$;

(ii) $g$ is positive on $[a,b]$ and even at $a$ and $b$ which means that $g$ is the restriction of a function $g$ such that $  g(a+t)=g(a-t), g(b+t)=g(b-t)$;

The proof is similar to the description of the metric in polar coordinates (see [Be] Lemma 9.114 and Theorem 9.125.)

 Note that it is  an easy exercise to prove that functions $f,g$ satisfying the ODE  characterizing $\AC$-metrics together with the initial conditions $(a),(b)$ extend to respectively odd and even (with respect to the points $a,b$ ) real analytic functions   around $a,b$.
 The fact that the conditions $(a),(b)$ of Theorem 1 imply (i) and (ii)  is based on the following elementary Lemma.

\medskip
{\bf Lemma 2.} {\it Let us assume that the function $Q$ is real analytic in a closed interval
$[x_0,x_1]$, $Q$ is positive in $(x_0,x_1)$ and $Q(x_0)=Q(x_1)=0,Q'(x_0)\ne 0,Q'(x_1)\ne 0$.
Let  a real function $\phi(t)$ satisfy an equation
$$\phi''=\frac12 Q'(\phi),\phi(0)=x_1,\phi'(0)=0,\tag 1.6$$
Then  $(\phi')^2=Q(\phi)$
 and  $\phi$ is a real analytic and periodic function $\phi\in C^{\omega}(\Bbb R)$ satisfying conditions $\phi(0)=x_0,\phi(l)=x_1$ where $l$ is the first positive real number such that $\phi(l)=x_1$ . Also $\phi$ is an even function at $t=0$ and $t=l$.}
\medskip
{\it Proof.} The first part of Lemma is well known (see [B],p.445, Lemma 16.37). If $\phi$ satisfies equation (1.6) then $(\phi')^2=Q_0(\phi)$ where $Q_0$ is a primitive function of $Q'$.
From the boundary conditions $Q_0=Q$.
 It is enough to show that $\phi$ is even at $t=0$ (for $t=l$ the proof is similar) which is equivalent to $\phi^{(2k-1)}(0)=0$ for every $k\in \Bbb N$. We prove  by induction with respect to $k$ that for every $F\in C^{\omega}(\Bbb R)$ the function $\psi(t)=F(\phi(t))$ is even at $t=0$. If $k=1$ then $\psi'(0)=F'(\phi(0))\phi'(0)=0$.  Let us assume that the result holds for $2k-1$. Let us write $\frac12Q'=Q_1$.
We get
$$\gather
\psi^{(2k+1)}=(\psi'')^{(2k-1)}=(F''(\phi)(\phi')^2+F'(\phi)\phi'')^{(2k-1)}\\=(F''(\phi)Q(\phi)+F'(\phi)Q_1(\phi))^{(2k-1)}=\\ \sum \frac{(2k-1)!}{p!(2k-1-p)!}((F''(\phi))^{(p)}Q(\phi)^{(2k-1-p)}+(F'(\phi))^{(p)}(Q_1(\phi))^{(2k-1-p)}).
\endgather$$
Thus $\psi^{(2k+1}(0)=0$. It follows that $F(\phi)$ is even at $0$ for every $F$. If $F=id$
then we get that $\phi$ is even at $0$.$\k$

Using this lemma one can show that in every case considered below in the paper function $g$
is even and $f$ is odd. Both functions are non-negative and $g$ is positive which follows from
the formula on $g$.

\bigskip

{\bf 2. Neutral bi-Hermitian Gray surfaces.}  In this section we shall construct bi-Hermitian neutral  metrics  on ruled surfaces $M_{k,g}$ of genus $g$.  The ruled surface  $(M_{k,g},g)$  is locally of co-homogeneity 1 with respect to the group of all local isometries of $(M_{k,g},g)$ and an open, dense submanifold $(U_{k,g},g)\subset (M_{k,g},g)$  is  isometric to the manifold
$ (a,b)\times P_k$
where $(P_k,g_k)$ is a 3-dimensional semi-Riemannian $\A$-manifold (a principal circle bundle $p:P_k\rightarrow \Si_g$ ) over a Riemannian surface $(\Si_g,g_{can})$ of constant sectional curvature  $K\in\{-4,0,4\}$ with a metric
$$g_{f,g}=dt^2+f(t)^2\th^2-g(t)^2p^*g_{can},\tag 2.1$$
where $g_k=\th^2-p^*g_{can}$ and $\th$ is the connection form of $P_k$ such that $d\th=2\pi k\  p^*\om$, where the de Rham cohomology class $[\om]\in H^2(\Si_g,\Bbb R)$ defined by the form $\om$ is an integral class corresponding to the class $1\in H^2(\Si_g,\Bbb Z)=\Bbb Z$.  Let $\th^{\sharp}$
be a vector field dual to $\th$ with respect to $g_k$. Let us consider a local orthonormal frame $\{X,Y\}$ on $(\Si_g,g_{can})$ and let $X^h,Y^h$ be horizontal lifts of $X,Y$ with respect to $p:M_{k,g}\rightarrow \Si_g$ (i.e. $dt(X^h)=\th(X^h)=0$ and $p(X^h)=X$) and let $H=\frac{\partial}{\partial t}$.
Let us define two almost Hermitian structures $J,\bJ$ on $M$ as follows
$$JH=\frac1f\th^{\sharp},JX^h=Y^h,\ \bJ H=-\frac1f\th^{\sharp},\bJ X^h=Y^h.$$
\medskip
{\bf Proposition 3. } {\it Let $\De$ be a distribution spanned by the fields $\{\th^{\sharp},H\}$.
Then $\De$ is a non-degenerate totally geodesic foliation with respect to the metric $g_{f,g}$ where $g$ is a non constant function. Both structures $J,\bJ$ are Hermitian and $\De$ is contained in the nullity of $J$ and $\bJ$.
The distribution $\DE$ is umbilical with the mean curvature normal $\xi=-\n\ln g$.  Let $\lb, \mu$ be eigenvalues of the Ricci tensor $S$ of $g_{f,g}$ corresponding to eigendistributions $\De,\De^{\perp}$ respectively.  Then the following conditions are equivalent:

(a) There exists $D\in \Bbb R$ such that $\lb-\mu =Dg^2$,

(b) There exist $C,D\in \Bbb R$ such that $\mu =Dg^2-C$,

(c) $\lb-2\mu$ is constant,

(d) $(U_{k,g},g_{f,g})$ is a neutral bi-Hermitian Gray surface.}
\medskip
{\it Proof.  } The first assertion of  Proposition 3 is a consequence similar to Proposition 3 in [J-2]. Note that $\n\lb=H\lb H,\n\mu=H\mu H$. Consequently $tr_g\n S=\frac12\n\tau=(H\lb+H\mu)H$. On the other hand one can easily check that $tr_g\n S=2(\mu-\lb)\xi+H\lb H$. Thus

$$\n\mu=2(\lb-\mu)\n\ln g.$$

Now we prove that (a) $\Rightarrow$ (b). If (a) holds then $\n \mu=2Dg^2\frac{\n g}g=D\n g^2$. Thus $\n(\mu-Dg^2)=0$ which implies (b).

(b)$\Rightarrow$ (a). We have
$$-\frac{\n g}g2(\mu-\lb)=\n\mu=2Dg\n g,$$
and consequently $\n g(\frac{(\mu-\lb)}g-Dg)=0$ which is equivalent to (a).

(a)$\Rightarrow $(c). We have $\lb-\mu=Dg^2$ and consequently  $\n \mu=2Dg\n g=D\n g^2$. Thus $\n \lb=\n (\mu+Dg^2)=2D\n g^2$ and $\n\lb-2\n\mu=0$ which gives (c).

(c)$\Rightarrow$ (a). If $\n\lb=2\n\mu$ then $\n\lb=4(\lb-\mu)\frac{\n g}g$. Consequently $\n\lb-\n\mu=2(\lb-\mu)\frac{\n g}g$ and $\n\ln|\lb-\mu|=2\frac{\n g}g=2\n\ln g$, which means that $\n\ln|\lb-\mu|g^{-2}=0$ or $\lb=\mu$ on the whole of $M$. In the last case we obtain an Einstein metric.  It follows that $\ln\frac{|\lb-\mu|}{g^2}=\ln D$ for some $D\in\Bbb R_+$ or $\lb-\mu=0$ (D=0), which is equivalent to  (a).

(d)$\Leftrightarrow$(c).  This equivalence follows from [J-3].
$\k$
\medskip
{\bf Theorem  2. } {\it On any  ruled  surface $M_{k,g}$ of genus $g$ with $k>0$  there exist a  one-parameter family of neutral bi-Hermitian $\AC$-metrics  $\{g_{x}:x\in(0,\e_s)\}$, where $\e_s>0$ depends only on $g$ and $k$,  which consists of neutral bi-Hermitian Gray metrics on $M_{k,g}$.}

\medskip
{\it Proof.}
 Note that for the first Chern class $c_1(\Si_g)\in H^2(\Si_g,\Bbb Z)$ of the complex curve $\Si_g$ we have the relation $c_1(\Si_g)=\chi\a$,  where $\a\in H^2(\Si_g,\Bbb Z)$ is an indivisible integral class and $\chi=2-2g$ is the Euler characteristic of $\Si_g$. Let us write $s=\frac{2k}{|\chi|}$ if $g\ne 1$ and $s=k$ if $g=1$. Then it is easy to show using O'Neill formulas for a semi-Riemannian submersion (see [ON],[B],[S]) that the manifold $(M_{k,g},g)$ with the metric $g$ given by $(2.1)$ has the Ricci tensor with the following eigenvalues :
$$\gather
\lb_0=-2\frac{g''}g-\frac{f''}f,\tag 2.2a\\
\lb_1= -\frac{f''}f-2\frac{f'g'}{fg}+2s^2\frac{f^2}{g^4},\tag 2.2b\\
\lb_2= -\frac{g''}g-\frac{f'g'}{fg}-(\frac{g'}g)^2-2s^2\frac{f^2}{g^4}-\frac{K}{g^2},\tag 2.2c
\endgather$$
where $\lb_0,\lb_1,$ correspond to eigenfields $T=\frac d{dt},\th^{\sharp}$ and $\lb_2$ corresponds to a two-dimensional eigendistribution orthogonal to $T$ and $\th^{\sharp}$.
Note that the equations we have got are practically the same as equations for Gray manifolds obtained in the Riemannian case. The only difference is that we now have $-K$ instead of $K$, so the cases $K>0$ and $K<0$ are now reversed, which give the different geometric meaning to the equations in Riemannian and semi-Riemannian neutral cases.  We present here the method of solving equations (2.2) for the convenience of the reader and completeness, although the calculations are just the same as in [J-4].
If $(M,g)\in\AC$ is a bi-Hermitian Gray surface then $\lb_0=\lb_1=\lb$, and if we  denote $\mu=\lb_2$, Proposition 3 b implies  an equation $$\mu=Dg^2-C\tag2.3$$
 for some $D,C\in \Bbb R$. Since $\lb_0=\lb_1$ we get
$$f=\pm\frac{gg'}{\sqrt{s^2+Ag^2}}.\tag 2.4$$
 Using a homothety of the metric we can assume that  $A\in\{-1,0,1\}$. In the case $A= 0$ we get a neutral  K\"ahler metric and these metrics on compact complex surfaces we shall describe in section 4. So we restrict our considerations to the case $A\in\{-1,1\}$.   Now we introduce a function $h$ such that $h^2=s^2+Ag^2$. Note that im$  h\subset (-s,s)$ if $A=-1$, and im$  h\subset (s,\infty)$  if $A=1$.
Then $g=\sqrt{|s^2-h^2|}$. Let us introduce a function  $z$ such that $h'=\sqrt{z(h)}$.  Note that
$$ f=h'\text{   and    } f'=\frac12z'(h).\tag 2.5$$
It follows that equation (2.3) is equivalent to
$$z'(h)-z(h)\frac{s^2+h^2}{h(s^2-h^2)}=\frac{4\e}h+\frac{D(s^2-h^2)^2}h-\frac{C(s^2-h^2)}h,\tag2.6$$
where $\e =-\text{sgn} KA\in\{-1,0,1\}$.  It follows that
\define \h{(\frac hs)}
$$\gather  z(h)=(1-\h^2)^{-1}(-4\e\h^{2}-\frac {Ds^4}5\h^6+(Ds^4-\frac{Cs^2}3)\h^4+\tag 2.7\\+(2Cs^2-3Ds^4)\h^2-4\e+Cs^2-Ds^4+\frac Es \frac hs).\endgather  $$

Let us denote again $C=Cs^2,D=Ds^4,E=\frac Es$ and let
$$z_0(t)=(1-t^2)^{-1}(-4\e(1+ t^{2})+D(-\frac15t^6+t^4-3t^2-1)+C(-\frac13
t^4+2t^2+1)+Et).
\tag 2.8$$
Write
$$P(t)=-4\e t^{2}-\frac{D}5t^6+(D-\frac{C}3)t^4+(2C-3D)t^2+Et-4\e+C-D.\tag 2.9$$
Then $z_0(t)=\frac{P(t)}{1-t^2} $. Note that  $z(h)=z_0\h$ and $z'(h)=\frac1sz_0'\h$. In view of Th. 1. we are looking  for real numbers $x,y\in \Bbb R, x>y$ such that
$$\gather
z_0(x)=0,  z'_0(x)=-2s,\tag 2.10a \\
z_0(y)=0,  z'_0(y)=2s,\tag 2.10b \endgather$$
and $z_0(t)>0$ for $t\in(y,x)$.
  Note that equations (2.10a) are equivalent to
$$\gather
-4\e x^{2}-\frac{D}5x^6+(D-\frac{C}3)x^4+(2C-3D)x^2-4\e+C-D+Ex=0\tag 2.11a\\
-8\e x-\frac{6D}5x^5+4(D-\frac{C}3)x^3+2(2C-3D)x+E=-2s(1-x^2).\tag2.11b\endgather$$
Equations $(2.11)$ yield
$$\gather D=\frac{5(-3E-6s-24\e x+3Ex^2-12sx^2-8\e x^3+2s x^4)}{2(-1+x)x(1+x)(15+10x^2-x^4)},\tag 2.12a\\
C=\frac{3(5E+10s+80\e x+30 s x^2-10E x^2+5Ex^4-10sx^4-16\e x^5+2s x^6)}{2(-1+x)x(1+x)(-15-10x^2+x^4)}\tag 2.12b
\endgather$$
Solving in a similar way equations $(2.10b)$ one can see that there exists a function $z_0$ satisfying the equations (2.10) if
$$\gather
(x+y)(-4\e(-5 x+x^3+5 y+2x^2y-2xy^2-y^3)\tag 2.13\\+s(5+2x^3y+2xy^3+3y^2+3x^2+x^2y^2-16xy))=0,\endgather$$
where  $x>y$, $x,y\in (-1,1)$ in the case $A=-1$ and $x,y\in (1,\infty)$ in the case $A=1$.
Using standard methods one can check that in the case of the genus $g\le 1$ (i.e. if $K=
4$ or $K=0$)
the only solutions of (2.13) giving a positive function $z$ are $x=-y\in(0,1)$. In the case of genus $g\ge 2$ we have
$K=-4$ and apart from the solutions with $x=-y$ (see [J-3]) there is an  additional family of solutions with $\e=-1$ for $s\in (0,\eta)$ on $M_{k,g}$ where $\eta=\frac2{17}(15+4\sqrt{13})\sqrt{\frac13(13(10\sqrt{13}-35))})\simeq
 2.05318..$. One can check that there exist families of solutions with $-1<y<x<1$ for $k\in\Bbb N$ such that $k<(g-1)\eta$. In fact it is not difficult to check that if $g>1$ then  on the sets $F=\{(x,y):-1<y<x<1\}$ and $H=\{(x,y):1<y<x\}$ the function $G(x,y)=-4\e(-5 x+x^3+5 y+2x^2y-2xy^2-y^3)+s(5+2x^3y+2xy^3+3y^2+3x^2+x^2y^2-16xy)$ is negative and positive somewhere on the boundaries of both  of $F$
and $H$  for
$s<2$. In the case $g\le1$ $G$ is positive on both $F$ and $H$. We have
$$\gather
G_x(x,y)=-4\e(-5+3x^2+4xy-2y^2)+2s(3x^2y+y^2x+3x-8y+y^3),\\
G_y(x,y)=-4\e(5+2x^2-4xy-3y^2)+2s(x^3+3y+x^2y-8x+3x^2y).\endgather$$
The equations $G_x=0,G_y=0$ are equivalent to
$$\gather
-4\e(-5+3x^2+4xy-2y^2)+2s(3x^2y+y^2x+3x-8y+y^3)=0,\\
(x+y)(10\e(y-x)+s(-5+x^2+3xy+y^2))=0.\endgather$$
Consequently for $\e=-1$ and  $s\in[2,\eta)$ $G$ attains its negative minimum in $F$ on the line $x=-y$.
It follows that if $g\le 1$ then $x=-y\in(0,1),E=0$ and $\e=1$ or $\e=0$. These solutions are also valid for $g\ge 2$ but in this case we have also solutions with $x\ne -y$ both for $\e=-1$ and $s\in(0,\eta)$ (then $-1<y<x<1$) and for $\e=1$ and $s\in(0,2)$ (in this case there exist solutions with $1<y<x$). Now we give explicit formulas for the case $x=-y\in(0,1)$.
Then we get

$$\gather
P(t)\tag 2.14\\=\frac1{x(15-5x^2-11x^4+x^6)}((t^2-x^2)(s(-15+10x^2-3x^4+t^2(10+12x^2-6x^4)\\
+t^4(-3-6x^2+x^4))+4\e x(x^2(-5+x^2)-t^4(3+x^2)+t^2(5+2x^2+x^4)))).\endgather$$
Thus $$P(0)=\frac{-4\e x^4(x^2-5)+sx(15-10x^2+3x^4)}{15-5x^2-11x^4+x^6}$$ and, since $\lim_{x\rightarrow 0^+}\frac{P(0)}x=s>0$,  there exists $\e_s>0$ such that  $P(t)>0$ if $t\in(0,x)$  for all $x\in(0,\e_s)$.  In fact in the case $\e=-1$ the real number $\e_s$ is the first positive root of the polynomial $-4x^3(x^2-5)+s(-15+10x^2-3x^4)$ and $\e_s=1$ if $\e=1$.
Note also that in both cases $\e_s=1$ if $s\ge 2$.
  Now the function $z_0(t)=\frac1{(1-t^2)}P(t)$ is positive on $(-x,x),x\in(0,\e_s)$. If $x\in(0,\e_s)$ then there exists a solution $h :(-a,a)\rightarrow (-sx,sx)$, where
$$a=\lim_{t\rightarrow sx^-}\int^t_0\frac {dh}{\sqrt{z_0\h}},$$
of an equation
$$h'=\sqrt{z_0\h},$$
such that $h(-a)=-sx,h(a)=sx,h'(-a)=h'(a)=0,h''(-a)=1,h''(a)=-1$. It follows that functions $f=h',g=\sqrt{s^2-h^2}$ are smooth on $(-a,a)$ and satisfy the boundary conditions described in Th.2. Consequently the metric
$$g_x=dt^2+f(t)^2\th^2-g(t)^2p^*g_{can}, $$
on the manifold $(-a,a)\times P_k$ extends to the smooth metric on the compact ruled surface $M=P_k\times_{S^1}S^2$ which is a $2-$sphere bundle over Riemannian surface $\Si_g$. Note that $g(-a)=g(a)=s\sqrt{1-x^2}$. $\k$
\medskip
{\bf Theorem  3. } {\it On the  surfaces $M_{0,g}=\Bbb{CP}^1\times{\Sigma_g}$ where $g\ge 2$ there exists a one-parameter family $\{g_{\a}:\a>1\}$ of   bi-Hermitian $\AC$-metrics. The Ricci tensor $\rho=\rho_{\a}$ of $(\Sigma_0,g_{\a})$ is bi-Hermitian and has two eigenvalues, which are everywhere different.}
\bigskip

{\it Proof. } Analogous to [J-3].
Let us write $g'=\sqrt{P(g)}$. Then
$$g^2P''(g)-2P'(g)g-4P(g)+16+6Cg^2=0.\tag 2.15$$
Consequently
$$P(g)=\frac Ag+Bg^{4}+Cg^2+4,\tag2.16$$
where $A,B\in \Bbb{R}$ are arbitrary.
 Now let  $D=1$ and let us consider the equations (we are looking
for unknown real numbers $A,B,C$ and $(x,y)$ where   $0<y<x$)
$$\gather P(y)=0,P(x)=0,\tag 2.17\\
P'(y)=2,P'(x)=-2.\endgather$$
Then
$y=\frac{4(\a-1)(\a^2+3\a+1)}{\a(2\a^2+\a+2)}$, $x=\a y=\frac{4(\a-1)(\a^2+3\a+1)}{(2\a^2+\a+2)}$ where $\a>1$. Note that $x,y>0$ and $C>0,A,B<0$.
Let us consider an equation ($P=P_{\a},h=h_{\a}$ depend on the parameter $\a>1$)
$$\frac{d^2h}{dt^2}=\frac12P'(h),\ h'(0)=0,\  h(0)=y=\frac{4(\a-1)(\a^2+3\a+1)}{\a(2\a^2+\a+2)}.\tag 2.18$$
This equation is equivalent to (if $t\in D=\{t\ge 0:h'(t)\ge 0\}$)
$$\frac{dh}{dt}=\sqrt{P(h)}, \  \  h(0)=y=\frac{4(\a-1)(\a^2+3\a+1)}{\a(2\a^2+\a+2)}.\tag 2.19$$
It follows that $P=P_{\a}$, where $\a>1$, has exactly two positive roots $\{x,y\}$ and $P(t)>0$ if $t\in (x,y)$. Note that  equation $(2.18)$ admits a smooth periodic solution $h$ defined on the whole of $\Bbb{R}$ and such that $\text{im }h=[x,y]$.
 Now it is easy to check that $\lb = -10Bh^2-3C$ and $\mu=-5Bh^2-3C$. The tensor $\rho-\frac{\tau}3g$ is a Killing tensor with eigenvalues $C,5Bh^2+C$ corresponding to $\De,\DE$ respectively. Then we  obtain a one parameter family of bi-Hermitian $\AC$-metrics $\{g_{\a}:\a>1\}$ on $M_{0,g}$.
$\k$
\bigskip

{\bf 3. Einstein neutral bi-Hermitian surfaces. }   Let us note that the solutions $(M_{k,g},g_x)$ with $D(x)=0$ correspond to Einstein neutral surfaces.  Let us recall  that
$$ D=\frac{5(-6s-24\e x-12sx^2-8\e x^3+2s x^4)}{2(-1+x)x(1+x)(15+10x^2-x^4)}\tag 3.1$$

Consequently there exist Einstein metrics in the family of Gray metrics $(M_{k,g},g_x)$ if and only if
the equation
 $$-6s-24\e x-12sx^2-8\e x^3+2s x^4=0\tag 3.2$$
has a real root $x\in (0,1)$.   Now it is not difficult to check that in the case $g=0$ we have $\e=1$ and $s\in \Bbb N$ and that equation (3.2) does not have any solution in $(0,1)$. Similarly in the case $g=1$.  In the case $g>1$ we have $s=\frac k{g-1}>0$ and $\e=-1$. In that case equation $(3.2)$ has a real root in $(0,1)$ if and only if $s\in(0,2)$.  To show this, let us denote
$$Q(x)=-6s-24\e x-12sx^2-8\e x^3+2s x^4.$$
Then
$$Q''(x)=-24(s+2\e x-sx^2).$$
Consequently $Q''(x)<0$ in $(0,1)$ for $\e\ge0$. In the case where $\e=-1$, since $Q''(0)=-24s<0$ and $Q''(1)=48>0$, the equation $Q''(x)=0$ has exactly real  one root, say $\a_s$.
(In fact, $\a_s=s(\sqrt{1+s^2}+1)^{-1}$.) Then, we have $Q'(\a_s)(=16(2-(1+s^2)(\sqrt{1+s^2}+1)^{-1}>0$ if and only if $(0<)s<(\sqrt{3}\slash2)^{\frac12}(\sqrt{3}+1)$.
In particular, for $0<s<2(<(\sqrt{3}\slash2)^{\frac12}(\sqrt{3}+1)$, $Q'(x)\ge Q'(\a_s)>0$, and hence $Q(x)$ is monotone increasing, in $(0,1)$. Since $Q(0)=-6s<0$ and $Q(1)=32-16s$, it follows that $Q(x)$ has exactly one root in $(0,1)$ if $s\in(0,2)$. On the other hand, if $s/ge 2$, then, since $3+6x^2-x^4>0$ and $x^4<x^3$ in $(0,1)$, we see that

$$\gather
Q(x)=24x+8x^3-2s(3+6x^2-x^4)\le 24x+8x^3-4(3+6x^2-x^4)\\<24x+8x^3-12-24x^2+4x^3=12(x^3-2x^2+2x-1)=\\12(x^3-2x^2+2x-1)=12(x-1)(x^2-x+1)<0.\endgather$$

 Thus the  equation $Q(x)=0$ has exactly one real root $x\in(0,1)$ if and only if $s\in(0,2)$. Since $s=\frac k{g-1}$ it follows that the real root $x\in(0,1)$ exists if and only if $k\in\{1,2,3,..,2g-3\}$.
Note that for $D=0$ the polynomial $P$ is even and of degree $4$ with exactly one positive root.
Consequently we get an Einstein metric on $M_{k,g}$ for all
$k\in \{1,2,..,2g-3\}$.  Thus we obtain on the ruled surfaces $M_{k,g}$  of genus $g>1$ constructed above exactly $2g-3$ different Einstein neutral bi-Hermitian metrics, exactly one Einstein metric on every surface $M_{k,g}$ for $k\in\{1,2,...,2g-3\}$. Note that for $g=2$ we obtain only one metric corresponding to the Riemannian Einstein Bergery-Page metric on the first Hirzebruch surface $F_1$.  Consequently we have proved

\medskip
{\bf Theorem  4. } {\it On a ruled  surface $M_{k,g}$ of genus $g\ge 2$ for $k\in\{1,2,...,2g-3\}$ there exists an Einstein   bi-Hermitian non-K\"ahler neutral metric.}

\medskip

 Note that in the Riemannian case we have the Page metric which is only one Einstein co-homogeneity 1 Hermitian, non-K\"ahler metric. (see [P],[B],[LeB],[S]) .

\bigskip
{\bf 4.   Neutral K\"ahler Gray surfaces.} The solutions with $A=0$ analogously as in the Riemannian case give neutral K\"ahler Gray  surfaces. We shall prove
\medskip
{\bf Theorem  5. } {\it On a ruled  surface $M_{k,g}$ of genus $g\ge 2$ for $k\in\{1,2,...,2g-3\}$ there exists   bi-Hermitian K\"ahler neutral metric.}
\medskip
{\it Proof.} Since $A=0$ we get  $f=\frac{gg'}s$. Consequently from (2.3) we get an equation

$$-2\frac{g''}g-4(\frac{g'}g)^2-\frac{K}{g^2}=Dg^2+C.\tag 4.1$$

Let us write $g'=\sqrt{P(g)}$. Then we can rewrite (4.1) as
$$P'(g)+\frac4gP(g)+\frac Kg+Dg^3+Cg=0,\tag 4.2$$
with the boundary conditions
$$\frac12yP'(y)=s,\frac12xP'(x)=-s,\tag 4.3$$
where $g(a)=y<g(b)=x$, $P(x)=P(y)=0$ and $P(t)>0$ for $t\in(y,x)$.

Consequently $P(g)=-\frac D8g^4-\frac C6g^2+\frac E{g^4}-\frac{K}4.$ Now exactly as in [J-2] we see that the solution satisfying the boundary conditions exists if $K=-4$ and $(s+2)y^4+(s-2)x^4=0$.
This condition implies that $s<2,C=0$, $D>0$ and $y=\root{4}\of{\frac{2(2-s)}D}, x=\root{4}\of{\frac{2(2+s)}D}$ and $E=\frac{s^2-4}{2D}<0$. In fact it is not difficult to check that a positive solution satisfying the boundary condition exists if and only if $s\in(0,2)$ hence for $k<2(g-1)$ and that the metrics corresponding to the same $s$ with different $D$ are homothetic.$\k$

  Note that the situation in the case of neutral 4-manifolds is again quite different from Riemannian case where we have only one irreducible non-Einstein compact K\"ahler Gray surface  (see [A-C-G],[J-2]).

\bigskip
The author is very grateful to the referee for the valuable comments and pointing out a mistake in the first version of the paper.

\bigskip
\centerline{\bf References.}
[A]   V. Apostolov   {\it Generalized  Goldberg-Sachs theorems for pseudo-Riemannian four-manifolds}   J.  Geom. Phys.  27,(1998), 185-198.
\medskip
\cite{A-C-G} V. Apostolov, D. Calderbank and P. Gauduchon {\it The geometry of weakly selfdual K\"ahler surfaces} Compositio Math. 135,(2003),179-322.
\par
\medskip
\cite{B} B\'erard-Bergery, L.{\it Sur de nouvelles vari\'et\'es riemanniennes d'Einstein}, Institut \'Elie Cartan,6, Univ. Nancy, (1982), 1-60.
\par
\medskip
\cite{Be}  A. Besse, {\it Einstein manifolds}, Springer Verlag, 1987.
\par
\medskip
\cite{D}  M. Dunajski, {\it  Anti-self-dual four manifolds with a parallel real spinor}, Proc. R. Soc. London. A(2002) 458, 1205-1222
\par
\medskip
[G] A. Gray {\it Einstein-like manifolds which are not Einstein} Geom. Dedicata {\bf 7} (1978) 259-280.
\par
\medskip
\cite{J-1} W. Jelonek,  {\it Compact K\"ahler surfaces with harmonic anti-self-dual Weyl tensor}, Diff. Geom. and Appl.16 (2002),267-276
\par
\medskip
\cite{J-2} W. Jelonek,  {\it Extremal K\"ahler $\AC$-surfaces},  Bull. Belg. Math. Soc. - Simon Stevin, 9,(2002), 561-571.
\par
\medskip
\cite{J-3} W. Jelonek,  {\it Bi-Hermitian Gray surfaces },  (to appear in Pacific Math. J)
\par
\medskip
\cite{J-4} W. Jelonek,  {\it Bi-Hermitian Gray surfaces II},  Preprint 2003
\par
\medskip
\cite{K} H. Kamada, {\it Self-dual K\"ahler metrics of neutral signature on complex surfaces}   Tohoku Mathematical Publications, September 2002.
\par
\medskip
\cite{LeB} C. Lebrun {\it  Einstein metrics on complex surfaces}, in Geometry and Physics (Aarhus 1995) eds. J. Andersen, J. Dupont, H. Pedersen and A. Swann, Lect. notes in Pure Appl. Math., Marcel Dekker, 1996.
\par
\medskip
\cite{M} Y. Matsushita  {\it Fields of 2-planes and two kinds of almost complex structures on compact 4-dimensional manifolds} Math. Zeitschrift {\bf 207} (1991), 281-291.
par
\medskip
\cite{M-L} Y. Matsushita and P. Law {\it Hitchin-Thorpe Type Inequalities for Pseudo-Riemannian 4-Manifolds of Metric Signature $(+ + - -)$} Geom. Dedicata {\bf 87} (2001), 65-89.
\par
\medskip
\cite{M-S} B.Madsen, H. Pedersen, Y. Poon, A. Swann {\it Compact Einstein-Weyl manifolds with large symmetry group.} Duke Math. J. {\bf 88} (1997), 407-434.
\par
\medskip
\cite{ON} B. O'Neill {\it Semi-Riemannian Geometry}, Academic Press, 1983.
\par
\medskip
\cite{P} D. Page {\it A compact rotating gravitational instanton}, Phys. Lett.79 {\bf B} (1978), 235-238.
\par
\medskip
\cite{Pe} D. Petean {\it Indefinite Ka\"hler-Einstein Metrics on Compact Complex Surfaces},
Comm. Math. Phys. 189  (1997), 227-235.
\medskip
\cite{S} P. Sentenac ,{\it Construction d'une m\'etriques  d'Einstein  sur la somme de deux projectifs complexes de dimension 2}, G\'eom\'etrie riemannienne en dimension 4 ( S\'eminaire Arthur Besse 1978-1979) Cedic-Fernand Nathan, Paris (1981), pp. 292-307.
\noindent

Institute of Mathematics

Technical University of Cracow

 Warszawska 24

31-155 Krak/ow,POLAND.

E-mail address: wjelon\@usk.pk.edu.pl

\enddocument